\documentclass[12pt,noinfoline]{imsart}

\usepackage{amsmath}
\usepackage{amssymb}
\usepackage{amsthm}
\usepackage{amscd}
\usepackage{amsfonts}
\usepackage{graphicx}%
\usepackage{fancyhdr, geometry,natbib}
\usepackage{url, hyperref, color, verbatim}
\usepackage{mathtools}
\usepackage{tikz}
\usetikzlibrary{matrix}
\usepackage{url, hyperref, color, verbatim}
\usepackage[normalem]{ulem}
\usepackage[ruled,  lined,  commentsnumbered,norelsize]{algorithm2e}
\theoremstyle{plain} 
\newtheorem{thm}{Theorem}[section]

\newtheorem{cor}[thm]{Corollary}

\newtheorem{prop}[thm]{Proposition}
\theoremstyle{definition}
\newtheorem{defn}[thm]{Definition}

\newtheorem{rmk}[thm]{Remark}

\newcommand{\R}{\operatornamewithlimits{R}}

\usepackage{tikz}
\usepackage{stmaryrd}

\definecolor{babyblue}{rgb}{0.54, 0.81, 0.94}

\newcommand{\mH}{\mathcal H}

\newcommand{\mB}{\mathcal{B}}

\newcommand{\ideal}[1]{\left< {#1}\right>}  
\newcommand{\vi}{{\sf V}}  

\DeclareMathOperator{\init}{LT_{\prec}}
\newcommand{\Basis}{\mathrm{Basis}}
\newcommand{\BruteForce}{\textrm{BruteForce}}

\def\jdlqed{\vbox{\hrule \hbox{\vrule\hbox to
5pt{\vbox to 6pt{\vfil}\hfil}\vrule}\hrule}}

\begin{document}

\title{The Spark Randomizer: a learned randomized framework for computing Gr\"obner bases} 
\author{Shahrzad Jamshidi and Sonja Petrovi\'c}
\thanks{SJ is at Lake Forest College. SP is at Illinois Institute of Technology. Supported by DOE award \#1010629 and the Simons Foundation Collaboration Travel Gift \#854770.}

\begin{abstract}
We define a violator operator which captures the definition of a minimal Gr\"obner basis of an ideal. 
This construction places the problem of computing a Gr\"obner basis within the framework of violator spaces, introduced in 2008 by G{\"a}rtner, Matou{\v{s}}ek, R{\"u}st, and {\v{S}}kovro{\v{n}} in a different context. The key aspect which we use is their successful utilization of a Clarkson-style fast sampling algorithm from geometric optimization. 
Using the output of a machine learning algorithm, we combine the prediction of the size of a minimal Gr\"obner basis of an ideal with the Clarkson-style biased random sampling method to compute a Gr\"obner basis in expected runtime linear in the size of the violator space. 
\end{abstract}

\maketitle


\section{Introduction}

The problem of finding a good basis for a polynomial system is the non-linear analogue of the computation of a reduced row echelon form.   This problem is essential in statistics, optimization, and other fields of science and engineering that rely on solving systems of polynomial equations. 
A good basis is one that allows algorithms to run and return unique outputs, and a Gr\"obner basis is a canonical example; see  \cite{StWhatIsGB} for a high-level overview.   A Gr\"obner basis of an ideal guarantees a solution to the ideal membership problem and the problem of detecting system feasibility, through the use of the multivariate polynomial division algorithm. 

Since the mid 20th century, when Gr\"obner bases \citep{Buchberger65}   and standard  bases \citep{Hironaka64}     were introduced,  there has been a tremendous amount of activity in developing methods and software for solving polynomial systems. All mainstream computer algebra systems offer generic algorithms for computing Gr\"obner bases applicable to all kinds of input ideals. Despite this great success in the field, unfortunately, algebraic problems have bad worst-case complexity.
 \citeauthor{Buchberger1965Translated}'s groundbreaking algorithm is able to compute a Gr\"obner basis of {any} ideal, but its computational complexity is difficult to estimate due to a number of choices one can make during its implementation. 
 What is known, for example by \cite{dube},  is that the maximum of the degrees of polynomials in the reduced Gr\"obner basis is doubly exponential in the number of variables. 
 \cite{BARDET201549} offer a nice overview of the upper and lower bounds that provide the general doubly exponential complexity result. 
Along a different avenue,  \cite{DylanPhD} discovered a new selection strategy for processing polynomials in Buchberger's algorithm  using reinforcement learning. 

 In a few cases, specialized algorithms have been used to improve runtime, such as \cite{BePa08,BePa08a,clo}. 
For the special case of toric ideals, \cite{ShortRationalToric} provide a polynomial time algorithm for Gr\"obner bases of toric ideals, which uses Barvinok's short rational generating functions and it returns not a list of polynomials, but rather the rational generating function; see also \cite{latte} and \cite{4ti2}. 

What is remarkable is that all known approaches for computing Gr\"obner bases of \emph{general} ideals rely upon variants of Buchberger's approximately 60-year-old algorithm.
 Improvements on it, such as Faug\`{e}re's  famous F5 algorithm developed in \cite{faugereetal, faugere2014sparse}, leverage the fact that the computation is a generalization of Gaussian elimination; see  \cite{BARDET201549} for an overview of  its complexity. 
 As such, these methods construct nontrivial organizational techniques, for example, cleverly organizing monomials into large matrices, to judiciously perform Buchberger's key step: reduction of S-polynomials through multivariate division. 
  
Our approach is distinct from these methods because \emph{we rely on a cleverly biased random sampling method} rather than {cleverly organized versions} of Buchberger's algorithm. We take a departure from the standard algorithms based on S-polynomials and address the problem of computing Gr\"obner bases using violator spaces.  
As we will see, violator spaces allow for adaptation of \emph{randomized} algorithms to computational non-linear algebra.
Figure~\ref{Spark_Randomizer} provides a schematic view of the proposed computational framework. 
For the non-expert, the ideal notation, vocabulary, and the setup of violator spaces are explained in Section~\ref{sec:background} on page \pageref{notation}. 
\usetikzlibrary{shapes}
\begin{figure}[h]
\centering 
\resizebox{11cm}{!}{%
\begin{tikzpicture}
	\node[draw, inner sep = 7pt, fill = magenta!20]  (generating_set) at (-6, 0) {\large $\{f_1, \ldots, f_s \}$ and $\prec$};
	\node [rounded corners, draw, inner sep = 4pt, fill = babyblue!40, ellipse] (ML1) at (-1,1.5) {\large $ML \ 1$};
	\node [rounded corners, draw, inner sep = 4pt, fill = babyblue!40, ellipse] (ML2) at (-1,-1.5) {\large $ML \ 2$};
	\node [rounded corners, draw, inner sep = 30pt, fill = violet!40] (H) at (2.5,1.5) {\Large $\mH$};
	\node[draw, fill=yellow, inner sep = 7pt] (sample) at (2.5, -1.5) {\large $\{\hat g_1, \ldots, \hat g_k \}$};
	\node[draw, ellipse, draw, inner sep = 4pt, fill = babyblue!40] (check) at (6.25, -1.5) {\large Check};
	\node[draw, fill=green!20,  inner sep = 7pt] (groebner) at (10, -1.5) {\large  $\{g_1, \ldots, g_k \}$};
``	\draw [->, ultra thick] (generating_set.east) to [out=0,in=180] (H.215);
	\draw [->, ultra thick] (generating_set.east) to [out=0,in=180] (ML1.west);
	\draw [->, ultra thick] (generating_set.east) to [out=0,in=180] (ML2.west);
	\draw [->, ultra thick] (ML1.east)  to [out =0, in = 180] (H.west);
	\draw[->, ultra thick] (ML2.east) to[out=0, in=180] (sample.west); 
	\draw[->, ultra thick] (H.south) to[out=-90, in=90] (sample.north);
	\draw[->, ultra thick] (sample.east) to[out=0, in=180] (check.west);
	\draw[->, ultra thick] (check.north) to[out=90, in=0] (H.east);
	\draw[->, ultra thick] (check.east) to[out=0, in=180] (groebner.west);
	\node at (0.35,1.9){\large $m$};
	\node at (.35,-1.1){\large $k$};
	\node at (6.5,0){\large no};
	\node at (7.85, -1.1){\large yes};
	\node[red] (biased) at (3.9, -0.2){\large \it biased};
	\draw [->, thick, red] (biased.south)arc(-160:160:.75);
	\node at (-6, -1){\large input};
	\node at (10, -2.5){\large output};
	\node at (2.5, -2.5){\large sample};
\end{tikzpicture}
}
\caption{\footnotesize A schematic view of the Spark Randomizer framework for a learned/randomized Gr\"obner computation. Given a generating set, $\{f_1, \ldots, f_s\}$ and a monomial order $\prec$, we predict two values using machine learning: the total degree ($m$) and the cardinality of a minimal Gr\"obner basis ($k$). From the generating set and $m$, we  construct a set $\mH$, from which we take samples of size $k$. We continually resample (biased) until we construct a minimal Gr\"obner basis $\{g_1, \ldots, g_k\}$, which is checked using the violator space  primitive query.\\
Note that prediction of $k$ is critical, while that of $m$ can be bypassed if there is an alternative way to construct  $\mH$ for a specific problem instance. }\label{Spark_Randomizer}
\end{figure}
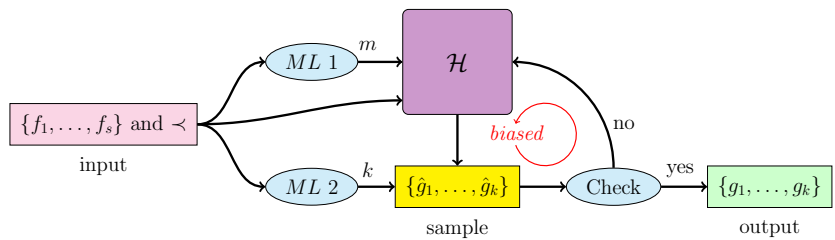

The input and output to the problem are standard: the input consists of a set of polynomials $f_1,\dots,f_s$ and a monomial order,  and the output is a minimal Gr\"obner basis of the ideal they generate. 
The two machine learning algorithms (``ML 1" and ``ML 2" in the Figure)  are used to predict the total degree ($m$) and the cardinality of the minimal Gr\"obner basis ($k$) of the ideal $\ideal{f_1,\dots,f_s}$; ML1 is optional. 
These algorithms guide the construction of the set $\mH$ and the size of each sample from it. 
We then use the framework of violator spaces and  Clarkson-style algorithms to sample from the universe $\mH$ and solve small  subproblems, embedded in an iterative biased sampling scheme. The small subproblems are solved using the violator space primitive query (``Check" in the Figure).  In the end,  the local information is used to make a global decision about the entire system. The expected runtime is linear in the number of input elements to the violator framework, $|\mH|$.

The Spark Randomizer blends randomization and learning with symbolic computation: 

\begin{thm}[The Spark Randomizer] \label{thm:main}
	 Let $I=\ideal{f_1,\dots,f_s}\subset K[x_1,\dots,x_n]$ be an ideal and fix a monomial ordering $\prec$. Let  $\mH\subset I$ be any finite set guaranteed to contain a minimal Gr\"obner basis of $I$. 
	 Let $k$ be the output of a (machine learning) algorithm which predicts the size  $\delta$ of a minimal $\prec$-Gr\"obner basis of $I$, with $k\geq \delta$. 

Then there exists a sampling algorithm that outputs a  Gr\"obner basis of $I$ in an expected number $O(k|\mH|+k^{O(k)})$ calls to the violator primitive query that solves a small monomial ideal membership problem. 

In particular, if  $k=\delta$,  the output is a minimal Gr\"obner basis. 
\end{thm}

This main result will follow as a consequence of the new violator space construction, Theorem~\ref{thm:grobner violator basis is a GB}. It is a combination of Corollary~\ref{cor: Clarkson for gb} and the outputs of machine learning algorithms, once we set up the appropriate violator space, provide an estimate of its combinatorial dimension $k$, and machine-learn  the sampling universe $\mH$ into existence. 
These three tasks are contents of Sections~\ref{sec:our violator}, \ref{sec:ML}, and \ref{section:H}, respectively.  
In particular, examples of worst-case universes $\mH$ appear on page \pageref{managable size of H}. 

Before providing the mathematical results, we summarize related work for further context. 

\subsection*{Why randomization?} 
In other areas of computational mathematics,  significant improvements in efficiency of deterministic algorithms have been obtained by algorithms that involve randomization. 
%
In 1992,  \citeauthor{sw-cblpr-92} identified special kinds of geometric optimization  problems that lend themselves to solution via repeated sampling of smaller subproblems: they called these {LP-type} problems. 
  A powerful sampling scheme,  devised in  \cite{c-lvali-95} for linear programming, works particularly well for geometric optimization problems in small number of variables.  Examples of applications include convex and linear programming, integer linear programming, the problem of computing the minimum-volume ball or ellipsoid enclosing a given point set in  $\R^n$, and the problem of finding the distance of two convex polytopes in $\R^n$.  
In 2008, G\"artner, Matou{\v{s}}ek, R\"ust and \v{S}kovro\v{n} \cite{ViolatorSpaces2008} invented {violator spaces} and showed they give a much more general framework to work with LP-type problems. In fact, violator spaces include all prior abstractions and were proven in \cite{SkovronPhDThesis} to be the most general framework in which Clarkson's sampling converges to a solution. It does so with expected runtime linear in the size of the input, as explained in Section~\ref{sec:background}. 

In symbolic computation specifically, randomization can be used to improve algorithm performance as illustrated, for example, in  \cite{AlgebraicViolators, primality, spielmanteng1, BSKW}.  

\subsection*{Violator spaces in non-linear algebra} 
There exists a precedent for using violator spaces for problems involving polynomial ideals. 
Namely,  \cite{AlgebraicViolators}  presented two violator operators: for computing   minimal generating sets  and certificate sub-systems  for large polynomial systems.  Their proof-of-concept implementation in {\tt Macaulay2} demonstrated the usefulness on an infeasible system built from a standard family of examples well-known for its Gr\"obner complexity, bringing down the computation time from $>20$ hours to just seconds.  
The randomized algorithm's    expected runtime is linear in the number of input polynomials. 
 Violator spaces make an appearance  in problems that have a natural linearization and a sampling size given by a combinatorial Helly number of the problem. 
  While violator spaces and Clarkson's algorithm have already a huge range of applications, to our knowledge, \cite{AlgebraicViolators}   was the first time such biased sampling  was used in computational algebraic geometry. 

\subsection*{Machine learning for non-linear algebra}  
Machine learning can be used to predict various ingredients in  symbolic algebraic computations. For example,  \cite{MLMathStructures} offers an overview of recent research on machine-learning mathematical structures; see also \cite{DeepLearnForSymbolicMath,DeepLearnForMath}. The thesis \cite{Silverstein}  uses features in neural network training to select the best algorithm to perform a Hilbert series computation by predicting a best choice of a monomial for one of the crucial steps of the computation, namely predicting a best pivot rule for the given input.

\cite{DylanMike-LearningBuch} and  \cite{MPP22}  demonstrate how Buchberger's algorithm is well-suited  for a machine learning approach, which is another precedent into which we lean. 
Specifically,  Buchberger's algorithm computes a Gr\"obner basis using an iterative procedure based on multivariate long division. The runtime of each step of the algorithm is typically dominated by a series of polynomial additions, and the total number of these additions is a hardware-independent performance metric that is often used to evaluate and optimize various implementation choices. \citeauthor{MPP22}  predict, using just the generating set, the number of polynomial additions that will take place during one run of Buchberger's algorithm. Good predictions are useful for quickly estimating difficulty and understanding what features make Gr\"obner basis computation hard. 
They show that a multiple linear regression model built from a set of easy-to-compute ideal generator statistics  {can predict the number of polynomial} {additions} somewhat well, better than an uninformed model, and better than regression models built on some intuitive commutative algebra invariants that are more difficult to compute. They also {train a simple recursive neural network that outperforms these linear models}.  
These results can directly be used for value models in the reinforcement learning approach to optimize Buchberger's algorithm introduced by \citeauthor{DylanMike-LearningBuch} who, within the paradigm of using learning to improve algorithms that give the exact answer, use machine learning to discover new S-pair selection
strategies in Buchberger's algorithm which outperform state-of-the-art human-designed heuristics by 20\% to 40\%.

Although these learning results are not used as direct input to Spark Randomizer---they are about Buchberger's algorithm heuristics---they  support the following core belief: that many algorithms in symbolic computation that relate to solving some polynomial systems are well-suited for machine learning techniques. The machine learning methods we do use for our framework are detailed in  Section~\ref{sec:ML}.

\section{Notation and background} \label{sec:background}

\label{notation}
Denote by  $K$  a   field, for example  the reader may keep $K=\mathbb C$ in mind.  Let $f_1=0,\dots,f_s=0$ be a system of $s$ polynomial equations in $n$ variables with coefficients in $K$.  
We will denote by  $\ideal{f_1,\dots,f_s}\subset R=K[x_1,\dots,x_n]$ the ideal generated by these polynomials. If $F=\{f_1,\dots,f_s\}$ is a set of polynomials, the ideal $\ideal{f_1,\dots,f_s}$ will equivalently be denoted by $\ideal{F}$. 
 
 Fix a monomial order $\prec$ on $R$ and an ideal $I$.  
Denote by $\init(h)$ the initial  term  of $h\in R$, which is the $\prec$-largest monomial. The initial ideal of $I$ with respect to $\prec$ is the monomial ideal $\init(I):=\ideal{\init(h):h\in I}$. 
 For the Gr\"obner basis definitions, \cite{clo} is a standard textbook reference. 

\begin{defn}
 A finite subset $B\subset I$ is a \emph{Gr\"obner basis} of $I\subset R$ with respect to $\prec$ if the initial terms of $B$ generate the entire initial ideal: $$\init(B):= \ideal{\init(b):b\in B}=\init(I):=\ideal{\init(f):f\in I}.$$
 A Gr\"obner basis is \emph{minimal} if the monomial generators of the initial ideal are minimal generators.  
 \end{defn}
 A minimal Gr\"obner basis is not unique. Its size, on the other hand, is always the number of minimal generators, or the ($0$-th) total  Betti number, of the  initial ideal $\init(I)$. 
 
 \cite{StWhatIsGB} offers a high-level overview of Gr\"obner bases and the textbook \cite{clo} discusses various applications as well. Of the many applications we single out two that have generated tremendous interest in recent decades: discrete optimization \cite{GBandOptim,rekha} and statistics \cite{GBandStats, St98}.


\label{sec:known}
\paragraph{}
From a different part of mathematics, we now introduce the basics of violator spaces, which were defined in the following context. An abstract LP-type problem is a tuple $(H,w,W,\leq)$. 
 A constraint $h\in H$  {\it violates\/} a set $G\subseteq H$ of constraints if  $w(G\cup\{h\})>w(G)$. 
The classical example  is the smallest enclosing ball problem: 
 $H$ is a finite point set in $\mathbb R^d$ and $w(G)$ is the radius of the smallest ball that encloses all points of $G$. 
  A point $h$ violates a set $G$ if it lies outside of the smallest ball  enclosing $G$. 
The {\it violator mapping} of $(H,w,W,\leq)$ is defined by
$\vi(G)=\{h\in H\colon w(G\cup\{h\})>w(G)\}$. Thus, $\vi(G)$ is the set 
of all constraints violating $G$. 
Finally, a violator space is defined as follows.

\begin{defn}[\cite{ViolatorSpaces2008}]
\label{defn:ViolatorSpaces}
Let $\mH$  be a finite set and $\vi$  a mapping $2^\mH\to2^\mH$. 
A {\it violator space} is a pair $(\mH,\vi)$ such that the following two axioms hold: 
\\\begin{tabular}{ll}
{\it Consistency}: & $G\cap \vi(G)=\emptyset$ holds for all $G\subseteq \mH$, and\\
{\it Locality}: & $\vi(G)=\vi(F)$ holds for all $F\subseteq G\subseteq H$ such that
$G\cap \vi(F)=\emptyset$.\\
\end{tabular}
\end{defn}

A basis of a violator space is defined in analogy to a basis of a linear programming problem: a minimal set of constraints that defines a solution space. Every proper subset of the basis is violated by a basis element. 
Violator spaces are also equipped with a natural combinatorial invariant, namely, the size of a largest basis. 
\begin{defn}[{\cite[Definitions~7 and 19]{ViolatorSpaces2008}}] \label{defn:delta}
A set $B\subseteq \mH$  
is a  {\it basis} if  $B\cap \vi(F)\neq\emptyset$ holds for all proper subsets $F\subset B$. 
For $G\subseteq \mH$, a {basis of $G$}  is a minimal subset $B$ of $G$ with $\vi(B)=\vi(G)$. 
The size of a largest basis of a violator space  $(\mH, V)$ is called the {\it combinatorial dimension} $\delta=\delta(\mH, V)$ of $(\mH, V)$.
\end{defn}

\citeauthor{ViolatorSpaces2008} proved that knowing the violations $\vi(G)$ for all $G \subseteq \mH$  is enough to compute the largest bases. To do so one can  utilize Clarkson's randomized algorithm, which is built on the primitive operation used as a black box in all stages of the algorithm: 
 \begin{defn}\label{defn:primitive} 
 Given a violator space $(\mH,\vi)$,  some set $G\subsetneq \mH$,  and some element $h\in\mH\setminus G$,  the  \emph{primitive test}  decides whether $h\in\vi(G)$. 
 \end{defn}

The following key result, which we will use, concerns the complexity  of finding a basis.   The runtime is given in terms of the combinatorial dimension  $\delta(\mH, V)$ and the size of $H$.  

\begin{thm}\cite[Theorem~27]{ViolatorSpaces2008} \label{key tool violators}
Using  Clarkson's algorithms, a basis of $\mH$ in a violator space $(\mH,V)$ whose combinatorial dimension is $\delta$ can be found by answering the  primitive query an expected $O\left(\delta \left|\mH\right| + \delta^{O(\delta)}\right)$  times.  
\end{thm}

\begin{rmk}\label{rmk: many primitive queries outperform large computation} 
The key insight is that answering the primitive query many times for small subsets of the universe is orders of magnitude faster, on average, than answering it once for the entire universe, which would be the brute-force method to solving the problem. 
\end{rmk}

The sampling method in \cite{c-lvali-95} avoids a full brute-force approach of testing each subset of size $\delta$ to see if it is a basis. The sampling method is presented  in two stages, referred to as  Clarkson's first and  second algorithm. 

\begin{algorithm}
\label{alg:ClarksonsFirst}
\LinesNumbered
\DontPrintSemicolon
\SetAlgoLined
\SetKwInOut{Input}{input}
\SetKwInOut{Output}{output}
\Input{
$G\subseteq \mH$, 
 			$\delta$: combinatorial complexity of $H$ }
\Output{ $\mB$, a basis for $G$ 
}
\BlankLine
\uIf{$|G|\leq9\delta^2$}{
return $\Basis2(G)$\;}
\Else{
$W\leftarrow\emptyset$\;
\Repeat{$V=\emptyset$}{
	$R\leftarrow$ random subset of $G\setminus W$ with $\lfloor\delta\sqrt{|G|}\rfloor$ elements.\;
	$C\leftarrow{\Basis2}(W\cup R)$\;
	$V\leftarrow\{h\in G\setminus C \text{ s.t.}\ h\in \vi(C)\}$\;
	\If{$|V|\leq 2\sqrt{|G|}$}{
		$W\leftarrow W\cup V$}
	}}
return $C$.
\caption{Clarkson's first algorithm, see also \cite{AlgebraicViolators}.}
\end{algorithm}

Clarkson's first algorithm, in the first iteration, draws a small random sample $R \subset G$, calls the second algorithm to calculate the basis $C$ of $R$, and returns $C$ if it is already a basis for the larger subset $G$. If $C$ is not already a basis, but the elements of $G\setminus C$ violating $R$ are few, it adds those elements to a growing set of violators $W$, and repeats the process with $C$ being calculated as the basis of the set $W\cup R$ for a new randomly chosen small $R\subset G\setminus W$. The crucial point here is that $|R|$ is much smaller than $|G|$. 

Clarkson's second algorithm ($\Basis2$) iteratively picks a random small ($6\delta^2$ elements) subset $R$ of $G$, finds a basis $C$ for $R$ by exhaustively testing each possible subset ($\BruteForce$), taking advantage of the fact that the sample $R$ is very small, 
and then calculates the violators of $G\setminus C$. 
At each iteration, elements that appear in bases with small violator sets get a higher probability of being selected. 

\begin{algorithm}
\label{alg:ClarksonsFirst}
\LinesNumbered
\DontPrintSemicolon
\SetAlgoLined
\SetKwInOut{Input}{input}
\SetKwInOut{Output}{output}
\Input{
$G\subseteq \mH$; $\delta$: combinatorial complexity of $H$. }
\Output{$\mB$: a basis of $G$}
\BlankLine
\uIf{$|G|\leq 6\delta^2$}{return \BruteForce$(G)$} 
\Else{\Repeat{$V=\emptyset$}{$R\leftarrow$ random subset of $G$ with $6\delta^2$ elements.\;
	$C\leftarrow\BruteForce(R)$\;
	$V\leftarrow\{h\in G\setminus C \text{ s.t.}\ h\in \vi(C)\}$\;
	\If{${\mathfrak m}(V)\leq {\mathfrak m}(G)/3\delta$}{
	\For{$h\in V$}{
		${\mathfrak m}(h)\leftarrow2{\mathfrak m}(h)$\;}}}}
return $C$.
\caption{Clarkson's second algorithm: $\Basis2(G)$, see  also \cite{AlgebraicViolators}.}
\end{algorithm}

It is critical to realize that the sampling process biased in a way that  some elements will be more likely to be chosen according to `how likely' they are to be elements in a basis. The algorithm does this by associating with every element $h$ of the set $G$  a multiplicity ${\mathfrak m}(h)$, and the multiplicity of a set is the sum of the multiplicities of its elements. The value $\mathfrak m(h)$ penalizes frequent appearance of $h$ as a violation of a given set. 

The entire process is repeated until a basis of $G$ is found, that is, until there are no violators of the selected (basis) set $C$. 


\section{Gr\"obner-violator framework}

\subsection{The violator operator that captures Gr\"obner bases} 
\label{sec:new} \label{sec:our violator}

To adapt this sampling scheme  to the problem of computing Gr\"obner bases, the key is to determine if there exists a violator space whose basis will be a Gr\"obner basis of an ideal. Given that violator  bases are minimal subsets that capture all violations of a given set, the following definition is intuitive. 
\begin{defn}
\label{defn:GbViolator}
Fix a monomial ordering $\prec$. 
Let $\mH\subset R$ be a finite set of polynomials and $S\subset\mH$. Define the  following  mapping  to capture the set of polynomials in the universe $\mH$ whose initial terms lie outside the initial ideal of S: 
 \[
 	\vi_{\prec}(S):=\{ h \in\mH : \init(S)\subsetneq\init(S\cup\{h\})\}.
 \]
 \end{defn}
 This is the natural violator space for the problem of computing a minimal Gr\"obner basis because the violations are those polynomials that grow the initial ideal within the universe $\mH$. 
 Proposition~\ref{prop:violator axioms} guarantees, for any  monomial order $\prec$ and any finite set $\mH\subset R$, that the pair $(\mH,\vi_\prec)$ is a violator space. 

  \begin{prop}\label{prop:violator axioms}
The operator $\vi_{\prec}(S)$  satisfies violator space axioms.  
 \end{prop}
 \begin{proof}
\emph{Consistency.}
Let $g\in G$. Then $\init(G)=\init(G\cup\{g\})$, thus $g\not\in \vi_\prec(G)$. 
On the other hand, let $g\in \vi_\prec(G)$. Then $\init(G)\subsetneq \init(G\cup\{g\})$. 
In particular, $\init(g)\not\in \init(G)$ and therefore $\init(g)\not\in \{\init(f):f\in G\}$, thus $g\not\in G$. 

\emph{Locality.} Suppose $F\subseteq G$ with $G\cap \vi_\prec(F)=\emptyset$. The aim is to prove that the violating sets $\vi_\prec(F)$ and $\vi_\prec(G)$ agree.  
Note that $F\subset G$ immediately implies $\vi_\prec(G)\subset\vi_\prec(F)$, so we need only prove the other inclusion. 

Given that no element of $G$ violates $F$,  
$	\init(F)=\init(F\cup\{g\})$
 holds for all $g\in G$, and in particular $\init(g)\in\init(F)$ so $\init(F)=\init(G)$. 

Suppose $h$ violates $F$, 
that is, $\init(F)\subsetneq \init(F\cup\{h\})$. Then the following sequence of inclusions holds:
\[
	\init(G)=  \init(\cup_{g\in G\setminus F} F\cup\{g\}) =  \init (F)  \subsetneq \init(F\cup\{h\}) \subseteq \init (G\cup\{h\}),
\]
and therefore $\init(G)\subsetneq \init(G\cup\{h\})$, in other words, $h$ also violates $G$. 
Therefore $\vi_\prec(F)\subset \vi_\prec(G)$, as required. 
\end{proof}

\begin{thm}\label{thm:grobner violator basis is a GB}

Fix a finite set of polynomials $\mH\subset R$  and a monomial order $\prec$.  
Then a largest basis of the violator space $(\mH,\vi_\prec)$ contains all the elements in the minimal Gr\"obner basis of $\ideal{\mH}$ which live in $\mH$. 
 
 In particular, for a fixed ideal $I\subset R$ and a finite subset $\mH\subset I$ that  is guaranteed to contain a $\prec$-Gr\"obner basis of $I$, a largest basis of $(\mH,\vi_\prec)$ is a minimal Gr\"obner basis of $I$. 
\end{thm}
\begin{proof}
A basis $B$ of a subset $G\subset\mH$ is an inclusion-minimal subset $B$ with the same violators as those of $G$. 
For the violator operator $\vi_\prec$, this translates to a set capturing the minimal generators of the initial ideal that live in $G$. One may think of any basis $B$ of $G\subset\mH$ as a partial minimal Gr\"obner basis of $I=\ideal{\mH}$: one that is built from elements in $G$ only. 

Since the universe $\mH$ is finite, there exists a largest basis $\mathcal B$ of the violator space $(\mH,\vi_\prec)$, namely, a set $\mathcal B$  such that $\mathcal B\cap \vi_\prec(F)\neq \emptyset$ for all proper subsets $F\subsetneq \mathcal B$. 
A largest basis $\mathcal B$ captures all violations of $\mH$ and consists of a set of polynomials in the universe $\mH$ whose initial terms minimally generate $\init(\mH)$.  

In the case when  $\mH\subset I$ contains a Gr\"obner basis of $I$, we have $I=\ideal{\mH}$, and therefore any largest basis of the violator space  is a minimal Gr\"obner basis of $I$. 

\end{proof} 

We can now apply Clarkson's algorithms  using the violator operator $V_\prec$ from Definition~\ref{defn:GbViolator}. 

\begin{cor}\label{cor: Clarkson for gb}
	 Let $I=\ideal{f_1,\dots,f_s}\subset K[x_1,\dots,x_n]$ be an ideal and fix a monomial ordering $\prec$. Let  $\mH\subset I$ be any finite set guaranteed to contain a $\prec$-Gr\"obner basis of $I$. Let $\delta$ be the combinatorial dimension of the violator space $(\mH,\vi_\prec)$. 
	 
Then there exists a sampling algorithm that outputs a minimal Gr\"obner basis of $I$ in an expected number $O(\delta|\mH|+\delta^{O(\delta)})$ calls to the primitive query that solves a small monomial ideal membership problem. 
\end{cor} 
Note that the  the combinatorial dimension  of the violator space $(\mH,\vi_\prec)$ is the size of any minimal $\prec$-Gr\"obner basis of the ideal $I$; see Definition~\ref{defn:delta}. 

The primitive query  for the operator $V_\prec$ decides whether an element $h\in\mH$ satisfies $h\in\vi_\prec(S)$  for a \emph{small} subset $S\subset \mH$; see Definition~\ref{defn:primitive}. 
  Running this primitive as a black box, as Clarkson does,  requires \emph{one} Gr\"obner basis computation to determine $\init(S)$, followed by the monomial ideal membership, which is a simple divisibility check. \label{primitive query}
We remind the reader of  Remark~\ref{rmk: many primitive queries outperform large computation},  which says that many primitive queries outperform one large computation of a Gr\"obner basis. This is a consequence of how the biased sampling algorithm is set up. 

Our results prove feasibility and correctness of the Spark Randomizer framework. 
What remains to make this framework \emph{practical} are the following: the estimate of, or upper bound on, the combinatorial dimension $\delta$, and the construction of the universe $\mH$ from which to sample. These are discussed in the next two sections; in particular,  we use  machine learning to predict the necessary combinatorial invariants.


\subsection{Learning the combinatorial dimension of the Gr\"obner violator space}\label{sec:ML}

There do exist upper bounds on the size of a reduced Gr\"obner basis appear in \cite{DylanMike-supplement}.  In the generic case  and for the grevlex  monomial ordering,  \cite{BARDET201549} provide a complexity result on the runtime of Faug\`{e}re's   F5 algorithm, and note that this is  also a bound on the number of polynomials for a reduced Gr\"obner basis, independently of the algorithm used. The sizes of a reduced and a minimal Gr\"obner bases agree and they equal the combinatorial dimension of the violator space $V_\prec$.  
In the context of Spark Randomizer, in reality, one only needs a `reasonable' upper bound on the combinatorial dimension $\delta$ for the algorithm to work correctly, where reasonableness is determined by an order of magnitude reduction from the size of the universe $\mH$. 
If the prediction is tight, the output is a minimal Gr\"obner basis; if it is an upper bound,  the output is a superset of it. 
 
 The question as to whether machine learning technology can reliably predict the cardinality and  degree of a reduced Gr\"obner basis given the input polynomials was initially tackled by \cite{PredictingGBcardinality2023}.  In that paper, the authors restrict their study to random binomial ideals with a fixed number of variables, generators, and total degrees, and they  train neural networks for the problem.
 They exploit a t-SNE visualization \cite{tSNE} to illustrate that it is difficult to classify the data for a lack of clustering and the absence of a decision boundary.
A t-SNE plot, short for t-distributed stochastic neighbor embedding, is a common statistical method for such a visualization. It assigns each datapoint a location on a 2- or 3-dimensional map; the data are then color-coded by the output, in the given case, size of the reduced graded reverse-lexicographic  Gr\"obner basis. 
The t-SNE algorithm relies on the Euclidean distance between these vectors and converting them to conditional probabilities in order to visualize data. This choice of distance is not an arbitrary choice precisely because  
  the $L^2$ norm, or mean squared error, is used in training neural networks via back propagation (i.e., gradient descent).
 
From the point of view of machine learning, the main difficulty lies in the fact that the Gr\"obner size prediction problems are very much \emph{unlike} traditional machine learning benchmark problems, in which the input is complex but the humans can quickly validate the prediction. 
\begin{figure}[h]
    \includegraphics[scale = 0.22]{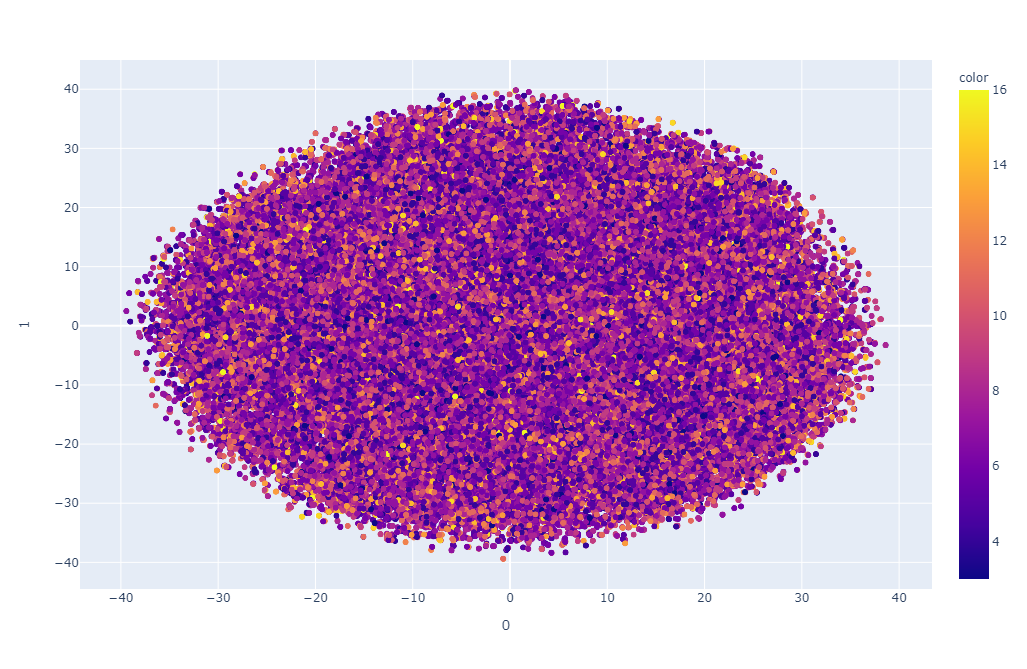}
 \caption{ A t-SNE plot for $99,722$ binomial ideals in $3$ variables, each generated by $5$ binomials (not necessarily homogeneous) with total degree $7$. Output color represents the size of the reduced graded reverse-lexicographic Gr\"ober basis.}
				  \label{fig:tSNEinhom}
\end{figure}
As an illustration, consider a data set of $99,722$ binomial ideals in $3$ variables, each generated by $5$ homogenous binomials with total degree $7$. For a human to `quickly validate the prediction' one either needs good folklore knowledge of the problem, which does not exist for a random ideal, or one needs a way to visualize this high-dimensional data in two or three dimensions.  
\citeauthor{PredictingGBcardinality2023} represent each ideal through the natural exponent vector representation  of the generators. 
The t-SNE plot  of this data set in Figure~\ref{fig:tSNEinhom}  demonstrates a lack of clustering and the absence of a decision boundary. This suggests a difficult classification problem.

  As there exist no benchmarks for this learning task, they compare the performance of their neural network again a simple linear regression model, and show  that the neural network model predicted the reduced Gr\"obner basis size better than multilinear regression from degree vector representations.  
\citeauthor{PredictingGBcardinality2023} generate and make available  a large data set that is able to capture sufficient variability in Gr\"obner complexity.
Neural networks were trained, on datasets of about one million binomial ideals, to predict the cardinality of a reduced Gr\"obner basis and the maximum total degree of its elements. 
 Simulations provided performance statistics on the order of $r^2 = 0.6$, which outperform both naive guess or multiple regression models which had  $r^2 = 0.180$. 
The $r$-squared statistic, or  the coefficient of determination, is a statistical measure of how close the data lie to the estimated regression line. For example, the baseline model which always predicts the mean will have an $r^2$ value of 0, and the model that exactly matches the observed values will have $r^2= 1$. Models that do worse than then baseline prediction will have a negative $r^2$ value; this often happens when the linear regression training and testing data sets are drawn from completely different distributions.

This demonstrates that the learning problem is highly nontrivial, but possible. 


\subsection{A universe for  the Gr\"obner violator and implications to complexity} \label{section:H}

The violator operator $V_\prec$ is a natural one for the Gr\"obner computation; however,  the universe from which we need to sample to find a basis does not follow the same kind of a setup as traditionally covered by Clarkson's algorithm. Namely, in problems of LP-type, such as finding the smallest enclosing ball for a finite set of given points in space, a basis of the violator space will be any  subset of points that defines the ball capturing the rest of the points. The key point is that  the universe $\mH$ is the set of given input points.  
In prior use of violators for computation with polynomials,  that the situation was  the same:  in \cite{AlgebraicViolators}, the large input system of polynomial equations is the universe from which one samples small subsets of polynomials to find a certificate of infeasibility. 
In contrast, the setup for the Spark Randomizer is different, in that the universe from which one draws small samples is not given to us by the input set of polynomials. Instead, we have to create it from the generating set.

Because the Spark Randomizer's runtime, on average, is linear in the size of the violator universe, it comes as no surprise that the complexity of the problem lives in the size of $\mH$.   Thus, crafting and pruning of $\mH$ is where the biggest gains of computational efficiency can be had.
To represent $\mH$ within the computer, we propose  a tensor representation of choices of monomials that can enter the ideal. 
  This can be stored recursively in tables---a kind of ad hoc database structure---relying on the specified order of the monomials.  A second critical practical question is pruning $\mH$. For example, if  one uses machine learning to estimate the maximum number of terms $\hat \gamma$ for any polynomial within the reduced Gr\"obner basis, and $\hat \gamma << \ell {d+n\choose n}$, then  $|\mH|$ can be reduced by deleting entries in the tensor. 

\paragraph{}
It is probably instructive to consider the worst case in several scenarios.  Recall that the number of monomials of degree exactly $d$ in $n$ variables is ${d+n\choose d}$, and the number of monomials up to degree $d$ is ${d+n\choose n}$. 
If $K=\mathbb F_\ell$ is a finite field with $\ell$ elements, then the total number of terms, $\gamma$, in each polynomial of degree at most $d$ in $R$ is at most  $\ell {d+n\choose n}$,  which grows as  $O({n+d\choose n}) = O((n+d)^{min(n,d)})$.  
In the homogeneous case, the upper bound on $|\mH|$ is $2^\gamma< 2^{\ell {d+n\choose d}}$, which for fixed $d$ smaller than $n$, this grows polynomially with the number variables $n$, on the order of $O(n^d)$. Of course, this upper bound is both  bad and  far from sharp, because the only ideal that will contain all these possible polynomials is $\ideal{1}$.  

For the special case of \emph{binomial} ideals, one could construct the universe $\mH$ by taking all pairs of  monomials up to degree $d$ to generate the universe of binomials, which can be done in time quadratic in the number of monomials, so $O({d+n\choose n}^2)$. This is, again, much larger than necessary as it does not take into account information from the input ideal; still, it results in a universe of theoretically manageable size:  $|\mH|=O(n^{2d}).$
\label{managable size of H}

\paragraph{}
For the sake of efficiency, crafting and storing all of $\mH$ before running the algorithm is probably not the most efficient idea. One should borrow insights from algebraic statistics and generate polynomials from the universe dynamically, rather than precomputing and storing the entire universe before running the sampling algorithm. Two such examples are  \cite{Dob2012} and \cite{GPS16} which, for a given matrix $A$, randomly generate binomials in the toric ideal $I_A$.  The former paper does so to generate minimal generating sets. The latter generates a superset of the Graver basis of  a $0/1$ matrix $A$ using the combinatorial structure of the hypergraph whose incidence matrix is $A$;  by construction, the output contains a valid $\mH$ for the Spark Randomizer, because the Graver basis contains all Gr\"obner bases of $I_A$ for every monomial order. In a shortly forthcoming joint paper, Bakenhus and Petrovi\'c  describe an efficient algorithm for generating large sets of binomials from any toric ideal so that the probability of generating every binomial in the Graver basis is nonzero. Its input is the matrix $A$ and the runtime is linear in the number of columns, which is the number of variables of the toric ideal.  
%
%
%

\section{Problems and discussion} 

The Spark Randomizer provides a feasible and correct framework for the use of randomization and machine learning to compute a Gr\"obner basis. Naturally, one can pose many questions to make the idea faster, more practical, or more flexible on specific families of examples. We outline a few such questions here. 

 \paragraph{}
The neural network model used in \cite{PredictingGBcardinality2023} predicted the reduced Grobner basis size better than multilinear regression. It performed especially well when the total maximum degree was fixed at 15. This is, in fact, exactly the outcome one would have hoped for—the learning problem is highly nontrivial, but possible. To a machine learning expert, both network architecture and algorithms used therein are  straightforward. 
On the other hand, the authors found no neural network model that predicted the maximum total degree well, nor one that could predict from summary statistics better than multiple linear regression,  which had a very low $r^2$ score (on the order of $0.1$ or negative). 
\paragraph{Problem 1.} Train a neural network model to predict maximum total degree of a minimal Gr\"obner basis element for a given ideal and a given monomial ordering.  

\paragraph{}
In all cases, there are likely more opportunities for improvement in tackling the learning problem by reconsidering some of the more complex hyperparameters in a neural network, such as the cost function. A cost function that is more meaningful to the problem could allow for significantly better convergence to an optimal outcome. 
The reader is referred to \cite[Section 3]{PredictingGBcardinality2023} for background on neural networks and an explanation of parameter choices. 

A slightly different point of view is to take a set of summary statistics from ideal generators as input to a learning algorithm, rather than all of the information, that is, explicit polynomials in full form. One could say there is a precedent for attempting such a summary: \citeauthor{MPP22} used a subset of features from a data set on binomial ideals which  did a reasonably good job in predicting the number of additions during one run of Buchberger's algorithm. \citeauthor{PredictingGBcardinality2023} failed to obtain similar results with \emph{the same} set of summary statistics; but this is not a surprise, because  performance metrics of Buchberger's algorithm need not correlate with complexity of Gr\"obner bases. Therefore, we ask: 

\paragraph{Problem 2.} Find a set of summary statistics, computable (in polynomial time) from the defining set of polynomials, that can be used to predict the maximum total degree of a minimal Gr\"obner basis element  for a given monomial ordering. Does such a set of summary statistics exist? 

\paragraph{} 
There is a possibility that there exist other violator spaces that can be useful to inform the Gr\"obner computation. For example, Buchberger proved that the only polynomials one ever needs in a Gr\"obner basis are constructed as S-polynomials. During the algorithm, their remainders on division by the growing basis set are added to the basis. Thus it might be tempting to construct a universe $\mH$ consisting of all possible S-polynomials one can construct up to  some degree $D$. This approach will neccesitate both  an estimate of $D$ and a different violator operator, such as : $\vi_{\prec}(S):=\{ h \in\mH : \mbox{ at least one term of } h \mbox{ is not in } \init(S)\}$. While the operator is consistent, we have not checked if the locality axiom holds.
\paragraph{Problem 3.} Are there other variants of the Gr\"obner violator operator which might lead to more efficient universes for sampling?

\paragraph{}  
Related to the above considerations, we ask for conditions under which the Spark Randomizer will outperform Buchberger's algorithm in practice.  One answer to the this part  is when the \emph{expected lineage} of S-pairs is very long relative to the complexity of the Gr\"obner basis. 
In \cite{tgb.m2}, we define and provide a {\tt Macaulay2} package that can compute such lineages. 
\begin{defn}[Definition 1.1 from \cite{tgb.m2}] 
Let $G$ be a  Gr\"obner basis of $I=(f_0,\dots,f_k)$.  
 A \emph{lineage} of a polynomial in $G$  is a natural number, or an ordered pair of lineages, tracing      its history in the given Gr\"obner basis computation using Buchberger's algorithm. It is defined recursively as follows: 
 \begin{itemize}
 \item  For the starting generating set,  $Lineage(f_i) = i$, 
 \item For any subsequently created S-polynomial $S(f,g)$, the lineage of its remainder $r$ on division is the pair 
 $Lineage(r)=(Lineage(f),Lineage(g))$. 
 \end{itemize}
 
\end{defn}

To illustrate, suppose $I=(x^2 - y, x^3 - z)\subset \mathbb Q[x,y,z]$ with graded reverse lexicographic order. Then $Lineage(x^2-y)=0$ and $Lineage(x^3-z)=1$.   
 Two additional elements are  added to create a (non-minimal) Gr\"obner basis: $xy + z$ and $y^2 -xz$, with lineages $(0,1)$ and $((0,1),0)$, respectively.  According to $Lineage(y^2-xz)$, this element is constructed from  $S(xy + z, x^2 - y)$. 
 Lineages are expressions of the starting generating set and thus dependent on the choice and order of its elements. 

The online documentation for our {\tt Macaulay2} package {\tt ThreadedGB.m2} offers an example of an S-polynomial with a long lineage but which becomes $1$, turning the reduced Gr\"obner basis into unity. 

\paragraph{Problem 4.} Predict the length of the longest lineage in one run of Buchberger's algorithm from the input set of polynomials, a monomial ordering, and the S-polynomial processing order that the algorithm will use. 

\paragraph{}
This problem, as stated, is a learning question that does not seem to be useful as direct input to the Spark Randomizer but, rather,  to inform the theoretical understanding of algorithm complexity and comparison. 
On the other hand,  recall that the violator space primitive query requires a small Gr\"bner computation, as described on page \pageref{primitive query}. With the lineage information learned, this computation can be organized more cleverly to minimize the expected lineage length and further speed up the Spark Randomizer by making the black-box primitive more efficient.

\paragraph{}
In general, understanding how learning can  further improve the algebraic versions of Clarkson’s algorithm  is an open problem wherein lies a promise of  much of the computational gain.

\bibliographystyle{apalike}
\bibliography{randomized-ideals,CompMathWhitepaper,MLandAlg,AlgStatAndNtwksAndMB}

\begin{thebibliography}{}

\bibitem[4ti2, 2018]{4ti2}
4ti2 (2018).
\newblock 4ti2-a software package for algebraic, geometric and combinatorial
  problems on linear spaces combinatorial problems on linear spaces.

\bibitem[Bardet et~al., 2015]{BARDET201549}
Bardet, M., Faug{\`e}re, J.-C., and Salvy, B. (2015).
\newblock On the complexity of the {F5} {G}r{\"o}bner basis algorithm.
\newblock {\em Journal of Symbolic Computation}, 70:49--70.

\bibitem[Beltr\'an and Pardo, 2008]{BePa08}
Beltr\'an, C. and Pardo, L. (2008).
\newblock On {S}male's 17th problem: a probabilistic positive solution.
\newblock {\em Foundations of Computational Mathematics}, 8(1):1--43.

\bibitem[Beltr\'an and Pardo, 2009]{BePa08a}
Beltr\'an, C. and Pardo, L. (2009.).
\newblock Smale's 17th problem: average polynomial time to compute affine and
  projective solutions.
\newblock {\em J. Amer. Math. Soc.}, 22(2):363--385.

\bibitem[Breiding et~al., 2018]{BSKW}
Breiding, P., Sturmfels, B., {Kali\v{s}nik-Verov\v{s}ek}, S., and Weinstein, M.
  (2018).
\newblock Learning algebraic varieties from samples.
\newblock {\em Revista Matematica Complutense}, 31:545--593.

\bibitem[Buchberger, 1965]{Buchberger65}
Buchberger, B. (1965).
\newblock {\em Ein Algorithmus zum Auffinden der Basiselemente des
  Restklassenringes nach einem nulldimensionalen Polynomideal}.
\newblock PhD thesis, Leopold-Franzens University.

\bibitem[Buchberger, 2006]{Buchberger1965Translated}
Buchberger, B. (2006).
\newblock Bruno {B}uchberger's {PhD} thesis 1965: An algorithm for finding the
  basis elements of the residue class ring of a zero dimensional polynomial
  ideal ({E}nglish translation).
\newblock {\em Journal of Symbolic Computation}, 41(3-4):475--511.

\bibitem[Clarkson, 1995]{c-lvali-95}
Clarkson, K.~L. (1995).
\newblock {L}as {V}egas algorithms for linear and integer programming.
\newblock {\em Journal of the ACM}, 42(2):488--499.

\bibitem[Cox et~al., 2007]{clo}
Cox, D., Little, J., and O'Shea, D. (2007).
\newblock {\em Ideals, Varieties, and Algorithms: An Introduction to
  Computational Algebraic Geometry and Commutative Algebra}.
\newblock Springer.

\bibitem[{De Loera} et~al., 2004]{ShortRationalToric}
{De Loera}, J., Haws, D., Hemmecke, R., Huggins, P., Sturmfels, B., and
  Yoshida, R. (2004).
\newblock Short rational functions for toric algebra and applications.
\newblock {\em Journal of Symbolic Computation}, 38(2):959--973.

\bibitem[{De Loera} et~al., 2013]{GBandOptim}
{De Loera}, J.~A., Hemmecke, R., and K\"oppe, M. (2013).
\newblock {\em Algebraic and Geometric Ideas in the Theory of Discrete
  Optimization}.
\newblock {MOS-SIAM} Series in Optimization. SIAM.

\bibitem[De~Loera et~al., 2004]{latte}
De~Loera, J.~A., Hemmecke, R., Yoshida, R., and Tauzer, J. (2004).
\newblock Effective lattice point enumeration in rational convex polytopes.
\newblock {\em Journal of Symbolic Computation}, 38(4):1273--1302.

\bibitem[{De Loera} et~al., 2016]{AlgebraicViolators}
{De Loera}, J.~A., Petrovi\'c, S., and Stasi, D. (2016).
\newblock Random sampling in computational algebra: Helly numbers and violator
  spaces.
\newblock {\em Journal of Symbolic Computation}, 77:1--15.

\bibitem[Dobra, 2012]{Dob2012}
Dobra, A. (2012).
\newblock Dynamic {M}arkov bases.
\newblock {\em Journal of Computational and Graphical Statistics}, pages
  496--517.

\bibitem[Dube, 1990]{dube}
Dube, T. (1990).
\newblock The structure of polynomial ideals and {G}r\"obner bases.
\newblock {\em SIAM Journal of Computing}, 19(4):750--773.

\bibitem[Faug\'ere et~al., 1993]{faugereetal}
Faug\'ere, J.-C., Gianni, P., Lazard, D., and Mora, T. (1993).
\newblock Efficient computation of zero dimensional {G}r\"obner bases by change
  of ordering.
\newblock {\em Journal of Symbolic Computation}, 16(4):329--344.

\bibitem[Faug\'ere et~al., 2014]{faugere2014sparse}
Faug\'ere, J.-C., Spaenlehauer, P.-J., and Svartz., J. (2014).
\newblock Sparse {G}r\"obner bases: the unmixed case.
\newblock arXiv preprint arXiv:1402.7205.

\bibitem[G{\"a}rtner et~al., 2008]{ViolatorSpaces2008}
G{\"a}rtner, B., Matou{\v{s}}ek, J., R{\"u}st, L., and {\v{S}}kovro{\v{n}}, P.
  (2008).
\newblock Violator spaces: structure and algorithms.
\newblock {\em Discrete Appl. Math.}, 156(11):2124--2141.

\bibitem[Gross et~al., 2015]{GPS16}
Gross, E., Petrovi\'c, S., and Stasi, D. (2015).
\newblock {G}oodness-of-fit for log-linear network models: {D}ynamic {M}arkov
  bases using hypergraphs.
\newblock {\em Annals of the Institute of Statistical Mathematics}.
\newblock DOI: 10.1007/s10463-016-0560-2.

\bibitem[He, 2022]{MLMathStructures}
He, Y.-H. (2022).
\newblock Machine-learning mathematical structures.
\newblock {\em International Journal of Data Science in the Mathematical
  Sciences}, 1(1):1--25.

\bibitem[Higham and Higham, 2018]{DeepLearnForMath}
Higham, C.~F. and Higham, D.~J. (2018).
\newblock Deep learning: An introduction for applied mathematicians.
\newblock {\em https://arxiv.org/abs/1801.05894}.

\bibitem[Hironaka, 1964]{Hironaka64}
Hironaka, H. (1964).
\newblock Resolution of singularities of an algebraic variety over a field of
  characteristic zero.
\newblock {\em Annals of Mathematics, Second Series}, 79(1):109--326.

\bibitem[Jamshidi et~al., 2023]{PredictingGBcardinality2023}
Jamshidi, S., Kang, E., and Petrovi\'c, S. (2023+).
\newblock Predicting the cardinality of a reduced gr\"obner basis.
\newblock {\em arXiv preprint arXiv:2302.05364}.

\bibitem[Jamshidi and Petrovi\'c, 2021]{tgb.m2}
Jamshidi, S. and Petrovi\'c, S. (2021).
\newblock Threaded {G}r\"obner bases: a {M}acaulay2 package.
\newblock {\em Journal of Software for Algebra and Geometry}, 11:123--127.

\bibitem[Lample and Charton, 2020]{DeepLearnForSymbolicMath}
Lample, G. and Charton, F. (2020).
\newblock Deep learning for symbolic mathematics.
\newblock In {\em International Conference on Learning Representations}.

\bibitem[Mojsilovi\'c et~al., 2023]{MPP22}
Mojsilovi\'c, J., Peifer, D., and Petrovi\'c, S. (2023).
\newblock Learning a performance metric of {B}uchberger's algorithm.
\newblock {\em Involve, a Journal of Mathematics}, 16(2):227--248.

\bibitem[Peifer, 2021]{DylanPhD}
Peifer, D. (2021).
\newblock {\em Reinforcement learning in {B}uchberger's algorithm}.
\newblock PhD thesis, Cornell University, Department of Mathematics.

\bibitem[Peifer et~al., 2020a]{DylanMike-LearningBuch}
Peifer, D., Stillman, M., and Halpern-Leistner, D. (2020a).
\newblock Learning selection strategies in {B}uchberger's algorithm.
\newblock In {\em Proceedings of the 37th International Conference on Machine
  Learning (ICML 2020)}.

\bibitem[Peifer et~al., 2020b]{DylanMike-supplement}
Peifer, D., Stillman, M., and Halpern-Leistner, D. (2020b).
\newblock Learning selection strategies in {B}uchberger's algorithm:
  Supplementary material.
\newblock {\em Proceedings of the 37th International Conference on Machine
  Learning (ICML 2020)}.

\bibitem[Robbiano, 2011]{GBandStats}
Robbiano, L. (2011).
\newblock {\em Gr\"obner Bases and applications}, volume 251 of {\em London
  Mathematical Society Lecture Note Series}.
\newblock Cambridge University Press.

\bibitem[Sharir and Welzl, 1992]{sw-cblpr-92}
Sharir, M. and Welzl, E. (1992).
\newblock A combinatorial bound for linear programming and related problems.
\newblock In {\em Proc. 9th Symposium on Theoretical Aspects of Computer
  Science (STACS)}, volume 577 of {\em Lecture Notes in Computer Science},
  pages 569--579. Springer-Verlag.

\bibitem[Silverstein, 2019]{Silverstein}
Silverstein, L. (2019).
\newblock {\em Probability and Machine Learning in Combinatorial Commutative
  Algebra}.
\newblock PhD thesis, University of California Davis.

\bibitem[{\v{S}}kovro{\v{n}}, 2007]{SkovronPhDThesis}
{\v{S}}kovro{\v{n}}, P. (2007).
\newblock {\em Abstract models of optimization problems}.
\newblock PhD thesis, Charles University, Prague.

\bibitem[Solovay and Strasse, 1977]{primality}
Solovay, R. and Strasse, V. (1977).
\newblock A fast {M}onte-{C}arlo test for primality.
\newblock {\em SIAM Journal of Computing}, 6(1).

\bibitem[Spielman and Teng, 2001]{spielmanteng1}
Spielman, D. and Teng, S. (2001).
\newblock Smoothed analysis of algorithms: why the simplex algorithm usually
  takes polynomial time.
\newblock {\em Journal of the ACM}.

\bibitem[Sturmfels, 1991]{St98}
Sturmfels, B. (1991).
\newblock {\em Gr\"{o}bner Bases and Convex Polytopes}, volume~8 of {\em
  University Lecture Series}.
\newblock American Mathematical Society.

\bibitem[Sturmfels, 2005]{StWhatIsGB}
Sturmfels, B. (2005).
\newblock What is... a {G}r\"obner basis?
\newblock {\em Notices of the American Mathematical Society},
  52(10):1199--1200.

\bibitem[Thomas, 1995]{rekha}
Thomas, R.~R. (1995).
\newblock A geometric buchberger algorithm for integer programming.
\newblock {\em Mathematics of Operations Research}, 20:864--884.

\bibitem[van~der Maaten and Hinton, 2008]{tSNE}
van~der Maaten, L. and Hinton, G. (2008).
\newblock Visualizing data using t-sne.
\newblock {\em Journal of Machine Learning Research}, 9:2579--2605.

\end{thebibliography}

 \end{document}